# Fractional Order Fuzzy Control of Nuclear Reactor Power with Thermal-Hydraulic Effects in the Presence of Random Network Induced Delay and Sensor Noise having Long Range Dependence

Saptarshi Das[a,b]*, Indranil Pan[a,c], and Shantanu Das[d]

a) Department of Power Engineering, Jadavpur University, Salt Lake Campus, LB-8, Sector 3, Kolkata-700098, India.
b) Communications, Signal Processing and Control Group, School of Electronics and Computer Science, University of Southampton, Southampton SO17 1BJ, United Kingdom.
c) Centre for Energy Studies, Indian Institute of Technology Delhi, Hauz Khas, New Delhi 110016, India.
d) Reactor Control Division, Bhabha Atomic Research Centre, Mumbai-400085, India.

Authors' Emails:
saptarshi@pe.jusl.ac.in, s.das@soton.ac.uk (S. Das*)
indranil.jj@student.iitd.ac.in, indranil@pe.jusl.ac.in (I. Pan)
shantanu@magnum.barc.gov.in (Sh. Das)

**Abstract:**
Nonlinear state space modeling of a nuclear reactor has been done for the purpose of controlling its global power in load following mode. The nonlinear state space model has been linearized at different percentage of reactor powers and a novel fractional order (FO) fuzzy proportional integral derivative (PID) controller is designed using real coded Genetic Algorithm (GA) to control the reactor power level at various operating conditions. The effectiveness of using the fuzzy FOPID controller over conventional fuzzy PID controllers has been shown with numerical simulations. The controllers tuned with the highest power models are shown to work well at other operating conditions as well; over the lowest power model based design and hence are robust with respect to the changes in nuclear reactor operating power levels. This paper also analyzes the degradation of nuclear reactor power signal due to network induced random delays in shared communication network and due to sensor noise while being fed-back to the Reactor Regulating System (RRS). The effect of long range dependence (LRD) which is a practical consideration for the stochastic processes like network induced delay and sensor noise has been tackled by optimum tuning of FO fuzzy PID controllers using GA, while also taking the operating point shift into consideration.

**Keywords**: fractional order fuzzy PID controller; long-range dependence; network induced stochastic delay; nuclear reactor thermal-hydraulics; power level control; anti-persistent noise.

## 1. Introduction

Nuclear reactors are becoming increasingly popular as a supplement to thermal power plants in providing electrical power to the national grid. As such nuclear reactor control is of prime essence due to the safety issues and command following operations.

The point kinetic dynamical model of a nuclear reactor is inherently governed by few sets of nonlinear differential equations. Its linearized transfer function models are marginally stable with the magnitudes of the eigen-values differing widely [1]-[2]. It is well known that few time constants in an open loop model of a nuclear reactor are very small due to the effect of the prompt neutron jump and few are quite large due to the delayed neutron jump [1]. Hence at different time instants, different time constants are predominant which makes the control of such kind of process very difficult. In [2], [3] spatial control of such a nuclear reactor is done by dividing the reactor into different zones and designing controller for each zone. However in the present case, the focus is on power level tracking. Hence an input-output model is derived first and then PID type controllers are used to remove any possible steady state offset. In Yun *et al.* [4], model predictive control is used for load following operation of a nuclear power plant at various operating points. However, in [4] the command-following power tracking responses have very sharp jumps which is detrimental for the system as it may cause thermal shock to the mechanical elements. Quantitative feedback theory (QFT) has been used in Torabi *et al.* [5] to design robust controllers to handle the nonlinearity and uncertainties in the reactor model and obtain favorable performance over a wide range of operating conditions. Fuzzy logic based intelligent system has been recently used in various nuclear power plant applications like use of neuro-fuzzy systems in power plant transient identification [6], fuzzy inference in nuclear reliability problems [7] etc. Several researchers have applied combination of various computational intelligence paradigms together for efficient control of the dynamics of nuclear reactors like in [8]-[11]. Among several other approaches, supervisory control [12], adaptive estimator based dynamic sliding mode control [13], gain scheduled dynamic sliding mode [14], robust nonlinear model predictive control [15], QFT [16] are gaining popularity in control problems related to nuclear power plants. In the present paper, an advanced fuzzy logic based fractional PID type controller is used which has the capability of handling change in the dynamics of a nuclear reactor with operating power level.

In Saha *et al.* [17], a fractional order phase shaper augmented with an optimal PID controller is proposed which is capable of providing operating point invariant control of nuclear reactor power over simple PID controllers. Das *et al.* [18] developed fractional dynamical models for different regimes of reactor operation and designed a robust FOPID controller to handle the shift in operating conditions. Atkin and Altin [19] have shown that the fuzzy logic based controllers work well for reactor power control even in the presence of noise and process parameter variation. In Liu *et al.* [20], the gains of the PID controller are updated with fuzzy inferencing mechanism, using GA to provide efficient nuclear reactor power control. The present approach combines all the above philosophies of fuzzy logic and fractional order controllers to design a FO fuzzy PID type controller for control of reactor power using real coded GA, at wide range of operating points. Coban [21] used particle swarm optimization (PSO) algorithm to tune the membership functions (MFs) of the fuzzy logic controller for nuclear reactor power level control. Since tuning of the input-output scaling factors (SFs) of fuzzy controllers has more impact on the control performance than tuning the MFs [22], the present paper tunes the input-output SFs and the orders of the error derivative and integral for fixed fuzzy rule-base and membership function type using a popular global optimization technique i.e. the genetic algorithm.

The Point kinetics equation based nonlinear model of a 500 MW nuclear reactor is derived in this paper, taking into account the thermal effects on reactivity. The nonlinear state space is linearized around the steady state operating conditions corresponding to the various operating powers of the nuclear reactor. The parameters of the FO fuzzy PID controllers are then tuned with these linearized models to obtain a trade-off design between good set point tracking and reduced controller effort. The fuzzy FOPID controllers tuned with the linearized models of the reactor at highest/lowest power shows good set-point tracking ability though the system's dc-gain changes due to change in the range of operation. Similar multiple model based control scheme for nuclear reactors has been proposed by Hu *et al.* [23] and the mathematical model for the guidance of control rod movement is proposed by Alireza and Shirazi [24].

Evaluation of the performance of nuclear reactor control systems over communication networks is being a big concern now-a-days because a small amount of stochastic delay is capable of destabilizing a well-tuned control loop as shown by Pan *et al.* [25]. Das *et al.* [26] have given experimental evidence of the impact of the presence of network in the control of nuclear reactors. Therefore, the possible random delays and the packet dropouts which occur during long distance transmission of signal from the sensors to the controller or from the controller to the actuator must be considered in the controller design phase itself, so as to prevent an unforeseen destabilization of reactor control loop, causing possible catastrophic failures. From the reactor physics point of view the dynamics of a nuclear reactor are mainly governed by two different types neutron groups viz. prompt neutron and delayed neutron [1]. At the start up conditions of nuclear reactor or at set-point changes, the random delays are highly detrimental especially for the prompt neutron jump and may cause rapid growth of global power and cause thermal shocks to its elements.

In networked control system (NCS) there is a strong possibility for the packets containing control signal to get delayed because of the shared network medium over which it is being transmitted. These stochastic delays are intrinsically different from the conventional process delays. Process delays are generally large and constant but the network induced delays are stochastically varying and have more adverse effects on the performance of a well-tuned control system [22], [27]. It has been shown in [28] that the network induced delay in a Local Area Network (LAN) exhibits self-similarity or non-Gaussian dynamics. The motivation of the present work lies in the fact that if the self-similarity of these processes or the delay dynamics associated with the network packets can be estimated properly then it is possible to minimize and compensate its deleterious nature. Bhambhani *et al.* [29] and Tejado *et al.* [30], [31] have suggested that fractional order dynamics of network induced delays can be best handled using fractional order controllers. But the fractional order noise and stochastic self-similar network delay are intrinsically different in nature [32]. Since the latter does not convey any extra energy in the control loop, the FO delay-dynamics in control applications needs to be extensively investigated. The present paper shows efficient method of optimization based fractional order fuzzy PID controller tuning to handle the shift in nuclear reactor operating power level as well as network induced stochastic delay and sensor noise which have fractional order dynamics.

In a nuclear power plant, control signals are generally passed through the dedicated channel from the reactor house to the distantly located control room. But with the advent of cheap communication and off the shelf hardware, Ethernet as a shared medium is getting importance to close real time control loops, in big complex plants. The advantage of using NCS is reduced wiring, modularity and flexibility over the existing technology. Thus online monitoring of reactor data and feedback of control signal can be done easily with LAN. But the major drawback is that congestion occurs in the data transfer processes which leads to delayed control signals resulting in poor control performance. So the data in the control loop suffering from stochastic delays due to the network congestion should not be fed-back directly to the RRS and some filtering/delay-compensation or control techniques should be used to reduce the stochastic nature of these signals. If these stochastic delays are not compensated, the RRS may malfunction and also may lead to tripping of the reactor. Often the sensor noise and network induced delay have the long memory effect and can thus be modeled using fractional order dynamics. Therefore, it is quite challenging to device a control algorithm which is capable of enforcing command-following reactor operation, considering the operating point shifting of reactor in the presence of fractal nature of the network delay dynamics and sensor noise. For this purpose, we studied the comparative performance analysis with conventional fuzzy PID and fractional order fuzzy PID controller in the reactor control loop. In recent literatures it has been seen that fractional order controllers perform well when the objective function is noisy and time-varying [25], [33]. This is due to the fact the classical PID controllers having derivative operators amplifies the noise, since at sharp edges the integer order derivative of error function becomes singular. But at those sharp edges the fractional derivative may exist [34], implying a relatively low amplification of the stochastic phenomena, introduced in the error signal due to random network induced delay or sensor noise, thus resulting in a better control performance. Also, it is shown in Pan *et al.* [22] that fuzzy logic based PID controllers due to having a rough or approximate reasoning in it to generate the control action, compared to simple mathematical operators like the proportional, integral and derivative operators, are capable of suppressing stochastic phenomena in control loops causing from random delay dynamics due to the presence of shared communication medium. Therefore, it is logical to numerically study the comparative performances of the conventional fuzzy PID controller and its fractional order counterpart, proposed by Das *et al.* [35] which are the hybridization of both the concepts, discussed above, for efficient nuclear reactor power level control problems with consideration of communication network delay and sensor noise having different dynamical characteristics.

The rest of the paper is organized as follows. Section 2 delineates the point kinetics based state space modeling and the consequent linearization of the nuclear reactor model. Section 3 discusses the structure and tuning of the fuzzy $PI^\lambda D^\mu$ type controllers with details of the fuzzy inferencing mechanisms. Simulation studies for the fuzzy FOPID controller with change in reactor power are shown in Section 4. Section 5 introduces basics of network induced delay and sensor noise with long range dependence. Section 6 proposes stochastic optimization based improved fuzzy FOPID controller tuning technique to suppress the unwanted randomness in RRS loop. The paper ends in Section 7 with the conclusions followed by the references.

**2. State-space modeling of the nuclear reactor**

The reactor model is developed using the point kinetic equations with six groups of delayed neutrons. The reactivity feedback due to the fuel temperature and coolant temperature are also taken into account in the reactor model development. In present modeling, the spatial effects in the reactor dynamics are neglected unlike [2] i.e. a single node or point reactor is considered.

*2.1. Neutron kinetic model of a point reactor*

$$\frac{dn}{dt} = \frac{\rho - \beta}{\Lambda} n + \sum_{i=1}^{G} \lambda_i c_i \tag{1}$$

$$\frac{dc_i}{dt} = \frac{\beta_i}{\Lambda} n - \lambda_i c_i, \quad i = 1, 2, \cdots, G \tag{2}$$

where, $n$ : neutron density ($m^{-3}$),
$\rho$ : total reactivity,
$\beta = \sum_{i=1}^{G} \beta_i$ : delayed neutron fraction,
$\Lambda$ : prompt neutron lifetime ($s$),
$G$ : groups of delayed neutron,
$\lambda_i$ : delayed neutron decay constant for the i$^{th}$ delayed neutron group ($s^{-1}$),
$c_i$ : precursor concentration for the i$^{th}$ precursor group ($m^{-3}$).

*2.2. Fuel temperature model*

The energy balance equation of the fuel pellet yields (3) which implies that the rise in temperature of the pellet is equal to the difference of the generated power due to nuclear fission and the heat dissipated by convection from the fuel.

$$m_f c_{pf} \frac{dT_f}{dt} = P \cdot n - Ah \left( T_f - \frac{T_i + T_e}{2} \right) \tag{3}$$

The term '$n$' in equation (3) denotes neutron density, appearing in equations (1) and (2). Now let us define, $\mu_f = m_f c_{pf}$ (for the conductive part of the heat transfer equation) and $\Omega = Ah$ (for the convective part of the heat transfer equation). Therefore, (3) can be rewritten as:

$$\frac{dT_f}{dt} = \frac{P}{\mu_f} n - \frac{\Omega}{\mu_f} T_f + \frac{\Omega}{2\mu_f} T_i + \frac{\Omega}{2\mu_f} T_e \tag{4}$$

where, $m_f$ : mass of the fuel ($kg$),
$c_{pf}$ : specific heat of the fuel at constant pressure ($J/kg.°C$),
$P$ : reactor power ($W$), $A$ : active heat transfer area ($m^2$),
$h$ : fuel-to-coolant heat transfer co-efficient ($W/m^2.°C$),
$T_f$ : average fuel temperature ($°C$),
$T_i$ : inlet coolant temperature ($°C$),

$T_e$: outlet coolant temperature ($°C$),

$\mu_f$: total heat capacity of the fuel ($J/°C$).

## 2.3. Coolant temperature model

The fuel pellet is surrounded by a Zircalloy cladding and is inserted into the heavy water which is used to cool the pellet. The dominant mode of heat transfer is by convection. The energy balance equation of the coolant yields (5) which implies that the rise in temperature of the coolant is equal to the difference of the convective heat transfer from the fuel pellet and the change in heat content of the coolant between the inlet and the exit.

$$m_c c_{pc} \frac{dT_c}{dt} = \Omega\left(T_f - \frac{T_i + T_e}{2}\right) - \left(w_c c_{pc} T_e - w_c c_{pc} T_i\right) \tag{5}$$

Let us define, $\mu_c = m_c c_{pc}$ and $M_c = w_c c_{pc}$. Therefore, equation (5) can be rewritten as (6):

$$\frac{dT_c}{dt} = \frac{\Omega}{\mu_c} T_f - \left(\frac{2M_c + \Omega}{2\mu_c}\right) T_e + \left(\frac{2M_c - \Omega}{2\mu_c}\right) T_i \tag{6}$$

where, $m_c$: mass of the coolant in the core ($kg$),

$c_{pc}$: specific heat of the coolant at constant pressure ($J/kg.°C$),

$w_c$: mass flow rate of the coolant ($kg/s$),

$\mu_c$: total heat capacity of the coolant ($J/°C$),

$T_c$: average temperature of the coolant ($T_c = \frac{T_i + T_e}{2}$).

Therefore,

$$\frac{dT_c}{dt} = \frac{1}{2}\frac{dT_e}{dt} \quad [\text{since}, \ T_i = \text{constant}] \tag{7}$$

Replacing (7) in (6) produces (8) as:

$$\frac{dT_e}{dt} = \frac{2\Omega}{\mu_c} T_f - \left(\frac{2M_c + \Omega}{\mu_c}\right) T_e + \left(\frac{2M_c - \Omega}{\mu_c}\right) T_i \tag{8}$$

## 2.4. Total reactivity model

The reactivity change may be caused due to the effect of control rod movement or thermal effects on reactivity given by (9).

$$\begin{aligned}\rho &= \rho_{rod} + \alpha_f \left(T_f - T_{f0}\right) + \alpha_c \left(T_c - T_{c0}\right) \\ &= \rho_{rod} + \alpha_f \left(T_f - T_{f0}\right) + \frac{\alpha_c}{2}\left(T_i - T_{i0}\right) + \frac{\alpha_c}{2}\left(T_e - T_{e0}\right)\end{aligned} \tag{9}$$

where, $\rho_{rod}$: reactivity introduced in the core due to control rod movement,

$\alpha_f$: temperature co-efficient of reactivity for the fuel ($°C^{-1}$),

$T_{f0}$: initial average temperature of the fuel ($°C$),

$\alpha_c$ : temperature co-efficient of reactivity for the coolant ($°C^{-1}$),

$T_{c0}$ : initial average temperature of the coolant ($°C$),

$T_{i0}$ : initial inlet temperature of the coolant ($°C$),

$T_{e0}$ : initial outlet temperature of the coolant ($°C$).

*2.5. Linearization and state-space model development*

Considering small perturbation around the steady-state operating point, nonlinear system of equations (1), (2), (4), (8) can be linearized to produce the following four set of coupled linear differential equations (subscript "*r*" denotes relative values):

$$\begin{aligned}
\delta\dot{n}_r &= \frac{\rho_0 - \beta}{\Lambda} \cdot \delta n_r + \lambda \cdot \delta c_r = f_1 \\
\delta\dot{c}_r &= \frac{\beta}{\Lambda} \cdot \delta n_r - \lambda \cdot \delta c_r = f_2 \\
\delta\dot{T}_f &= \frac{P_0}{\mu_f} \cdot \delta n_r - \frac{\Omega}{\mu_f} \cdot \delta T_f + \frac{\Omega}{2\mu_f} \cdot \delta T_i + \frac{\Omega}{2\mu_f} \cdot \delta T_e = f_3 \\
\delta\dot{T}_e &= \frac{2\Omega}{\mu_c} \cdot \delta T_f + \left(\frac{2M_c - \Omega}{\mu_c}\right) \cdot \delta T_i - \left(\frac{2M_c + \Omega}{\mu_c}\right) \cdot \delta T_e = f_4
\end{aligned} \qquad (10)$$

At steady state, equation (2) reduces to

$$\dot{c}_{r0} = 0 \Rightarrow \frac{\beta}{\Lambda} \cdot n_{r0} - \lambda \cdot c_{r0} = 0 \Rightarrow c_{r0} = \frac{\beta}{\Lambda\lambda} \cdot n_{r0} \qquad (11)$$

At steady state equation (1) also reduces to

$$\dot{n}_{r0} = 0 \Rightarrow \frac{\rho_0 - \beta}{\Lambda} n_{r0} + \lambda c_{r0} = 0 \qquad (12)$$

In (11) and (12), the subscripts "0" correspond to the steady state values. Now, substitution of the steady state delayed neutron precursor concentration ($c_{r0}$) from (11) to (12) yields $\rho_0 = 0$. Since, inlet coolant temperature is generally kept constant, thus yielding $\delta T_i = 0$. Now, equation (10) can be converted to the conventional SISO state-space model for the reactor as follows:

$$\begin{aligned}
\dot{x}(t) &= Ax(t) + Bu(t) \\
y(t) &= Cx(t) + Du(t)
\end{aligned} \qquad (13)$$

Here, the input ($u$), output ($y$) and state variables ($x$) are defined as $u = \delta\rho_{rod}$ (change in reactivity of the core due to control rod movement), $y = \delta n_r$ (relative power output), $x = \begin{bmatrix} \delta n_r & \delta c_r & \delta T_f & \delta T_e \end{bmatrix}^T$. The system matrices are defined below in terms of the Jacobian of the four functions in (10):

$$A = \begin{bmatrix} \frac{\partial f_1}{\partial x_1} & \frac{\partial f_1}{\partial x_2} & \frac{\partial f_1}{\partial x_3} & \frac{\partial f_1}{\partial x_4} \\ \frac{\partial f_2}{\partial x_1} & \frac{\partial f_2}{\partial x_2} & \frac{\partial f_2}{\partial x_3} & \frac{\partial f_2}{\partial x_4} \\ \frac{\partial f_3}{\partial x_1} & \frac{\partial f_3}{\partial x_2} & \frac{\partial f_3}{\partial x_3} & \frac{\partial f_3}{\partial x_4} \\ \frac{\partial f_4}{\partial x_1} & \frac{\partial f_4}{\partial x_2} & \frac{\partial f_4}{\partial x_3} & \frac{\partial f_4}{\partial x_4} \end{bmatrix}; B = \begin{bmatrix} \frac{\partial f_1}{\partial u} \\ \frac{\partial f_2}{\partial u} \\ \frac{\partial f_3}{\partial u} \\ \frac{\partial f_4}{\partial u} \end{bmatrix}. \tag{14}$$

Therefore the reactor model can be written in the form (13) with the system matrices given by (15).

$$A = \begin{bmatrix} -\frac{\beta}{\Lambda} & \lambda & \frac{n_{r0}\alpha_f}{\Lambda} & \frac{n_{r0}\alpha_c}{2\Lambda} \\ \frac{\beta}{\Lambda} & -\lambda & 0 & 0 \\ \frac{P_0}{\mu_f} & 0 & -\frac{\Omega}{\mu_f} & \frac{\Omega}{2\mu_f} \\ 0 & 0 & \frac{2\Omega}{\mu_c} & -\frac{2M_c + \Omega}{\mu_c} \end{bmatrix}; B = \begin{bmatrix} \frac{n_{r0}}{\Lambda} \\ 0 \\ 0 \\ 0 \end{bmatrix}; \tag{15}$$

$$C = \begin{bmatrix} 1 & 0 & 0 & 0 \end{bmatrix}; D = \begin{bmatrix} 0 \end{bmatrix}.$$

The reactor parameters at various power levels are reported in Table 1 as studied in [4].

**Table 1: Reactor parameters at various operating power levels**

| Power (%) | $T_c$ | $n_{r0}$ | $\alpha_f$ | $\alpha_c$ | $\mu_f$ | $\mu_c$ | $\Omega$ | $M_c$ |
|---|---|---|---|---|---|---|---|---|
| 100 | 302 | 1.0 | $-2.9 \times 10^{-5}$ | $-6.3 \times 10^{-4}$ | $2.25 \times 10^7$ | $6.9 \times 10^7$ | $3.94 \times 10^6$ | $7.08 \times 10^7$ |
| 80 | 298.6 | 0.8 | $-3.2 \times 10^{-5}$ | $-5.59 \times 10^{-4}$ | $2.21 \times 10^7$ | $6.8 \times 10^7$ | $4.16 \times 10^6$ | $6.89 \times 10^7$ |
| 60 | 295 | 0.6 | $-3.3 \times 10^{-5}$ | $-5.56 \times 10^{-4}$ | $2.18 \times 10^7$ | $6.7 \times 10^7$ | $4.38 \times 10^6$ | $6.87 \times 10^7$ |
| 40 | 291.8 | 0.4 | $-3.5 \times 10^{-5}$ | $-5.22 \times 10^{-4}$ | $2.14 \times 10^7$ | $6.61 \times 10^7$ | $4.61 \times 10^6$ | $6.79 \times 10^7$ |
| 20 | 288.4 | 0.2 | $-3.8 \times 10^{-5}$ | $-4.86 \times 10^{-4}$ | $2.10 \times 10^7$ | $6.53 \times 10^7$ | $4.85 \times 10^6$ | $6.7 \times 10^7$ |

The full power of the nuclear reactor is considered to be $P_0 = 500 MW$. Other typical constants governing the reactor dynamics like the delayed neutron fraction, prompt neutron lifetime and delayed neutron decay constants have been considered similar to the study reported by Theler and Bonetto [36] as $\beta = 7.65 \times 10^{-3}$, $\Lambda = 1.76 \times 10^{-4}$, $\lambda = 7.59 \times 10^{-2}$. The state space models (15) around different reactor power can be converted to the corresponding transfer function models using the well-known relation (16).

$$G(s) = C(sI - A)^{-1} B + D \tag{16}$$

It is clear from the linear transfer function models of the nuclear reactor (17) that its dc gain widely varies with the shift in global reactor power. Also, the models have a pole very close to the origin indicating a marginally stable open loop dynamics with strong lead-time constants and high dc-gain. The co-existence of such complex dynamical behaviors makes the task of maintaining command following reactor control very difficult, especially taking the shift in operating point (reactor power) into consideration. Few previous studies on robust fractional order controller design for efficient nuclear reactor operation in different power levels have been done by Saha *et al.* [17] and Das *et al.* [18]. In the linearized models ($G$) of the nuclear reactor in equation (17) the subscripts denote different operating power levels.

$$G_{100} = \frac{5681.8182(s+2.114)(s+0.17)(s+0.0759)}{(s+43.52)(s+2.096)(s+0.1979)(s+0.01682)}$$

$$G_{80} = \frac{4545.4545(s+2.094)(s+0.1822)(s+0.0759)}{(s+43.52)(s+2.08)(s+0.2066)(s+0.0137)}$$

$$G_{60} = \frac{3409.0909(s+2.123)(s+0.1941)(s+0.0759)}{(s+43.52)(s+2.112)(s+0.2137)(s+0.01044)} \quad (17)$$

$$G_{40} = \frac{2272.7273(s+2.132)(s+0.2076)(s+0.0759)}{(s+43.53)(s+2.124)(s+0.2217)(s+0.007051)}$$

$$G_{20} = \frac{1136.3636(s+2.135)(s+0.2219)(s+0.0759)}{(s+43.53)(s+2.132)(s+0.2296)(s+0.003624)}$$

DC-gain of the linearized models of the nuclear reactor at different operating powers in (17) are 510.424, 513.7372, 519.9578, 528.2555, 529.1494 respectively which also justifies the frequency response shown in Figure 1, since stability of the reactor models decreases at lower powers due to increased dc-gain. This is evident from the gradual reduction in the open loop phase at low frequency regions in the Bode diagram and gradual shifting of the Nyquist curve towards the critical point (-1,0) for low reactor powers. Attempt for frequency domain design of robust FO controllers with varying level of reactor power can be found in Saha *et al.* [17] and Das *et al.* [18].

It is also to be noted that the reactor under consideration is of Pressurized Water (PWR) type. The thermal hydraulic parameters for a PWR have been adopted from [4]. The thermal hydraulic part of the state-space model is based on the transient heat balance for the fuel and coolant. For simple and small pool type research reactors generally the neutronics is only considered. In such cases, there is no thermal loop taking part in power controls, like power reactors of Canadian Deuterium Uranium (CANDU) type, PWR and Boiling Water Reactor (BWR) type etc. For small research reactors the idea is to govern the basic neutronics only, as done conventionally by ON-OFF controls of the control rods, looking at only the linear neutron channel signal [37], [38]. But for power reactors unlike research reactors have a loop of thermal-hydraulics and thus control inputs from these parts to affect the reactivity.

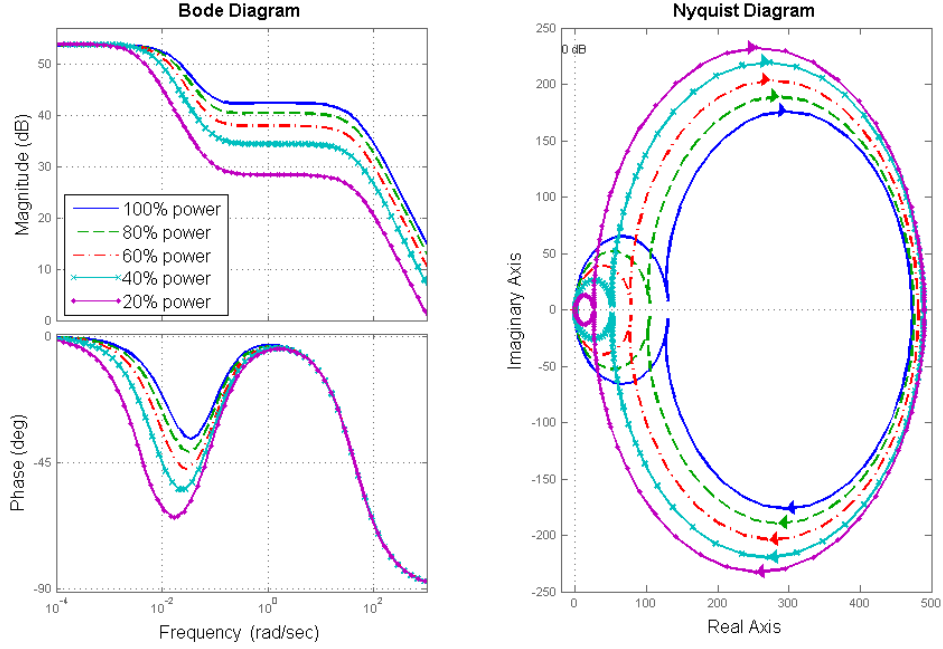

**Figure 1: Frequency response of the reactor, operating at different powers.**

## 3. Concept of Fuzzy FOPID controller and its optimum tuning
### *3.1. Basics of fractional order fuzzy PID controller*

Simple fuzzy PID controller and its global optimization based tuning to compensate for the network induced delays has been extensively studied in Pan *et al.* [22]. The idea has been extended by Das *et al.* [35] with its fractional order counterpart to handle nonlinear and open-loop unstable plants. For the present fuzzy FOPID controller $\{K_e, K_d\}$ are the input SFs and $\{K_{PI}, K_{PD}\}$ are the output SFs as shown in Figure 2. Here, the integer order rate of error in the conventional integer order fuzzy logic controller (FLC) input has been replaced by its FO counterpart ($\mu$). Also the FLC output is fractionally integrated with order ($\lambda$) which can be tuned to meet designer's specifications. The control law of this controller is given by (18).

$$u_{FLC\_FOPID}(t) = u_{FLC\_FOPI}(t) + u_{FLC\_FOPD}(t)$$
$$= K_{PI} \cdot \frac{d^{-\lambda} u_{FLC}(t)}{dt^{-\lambda}} + K_{PD} u_{FLC}(t) \tag{18}$$

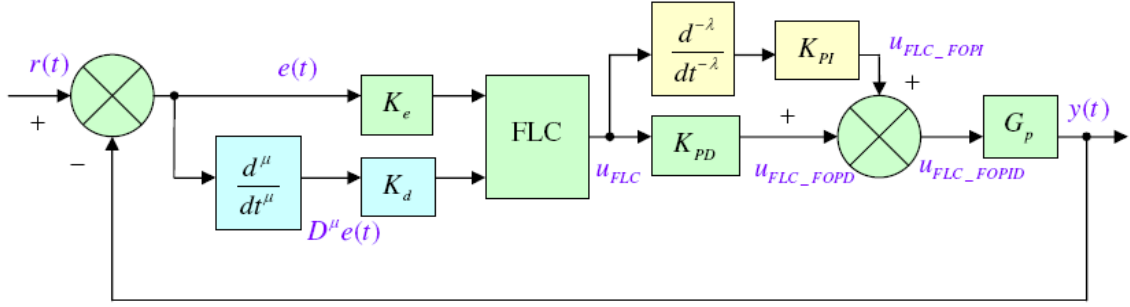

**Figure 2:** Schematic diagram of FO fuzzy PID controller.

For the typical integer order fuzzy PID controller in [22], [35] with product-sum inferencing, center of gravity defuzzification method and triangular membership function, the relation between the input and output variables can be expressed as

$$u(t) = \widetilde{A} + \widetilde{B}e(t) + \widetilde{C}\dot{e}(t) \tag{19}$$

where, parameters $\{\widetilde{A}, \widetilde{B}, \widetilde{C}\}$ are devised from the error and control signal. Therefore, it is expected that for the controller in Figure 2, the FLC output (19) will be a function of fractional rate of error instead of conventional integer order derivative of error signal. Now, using the well-known identity of fractional calculus or successive derivative of higher powers of a variable $\dfrac{d^n}{dt^n}t^m = \dfrac{\Gamma(m+1)}{\Gamma(m-n+1)}t^{m-n}$, equation (18) can be expressed as

$$\begin{aligned}u_{FLC\_FOPID}(t) &= K_{PI} \cdot \frac{d^{-\lambda}}{dt^{-\lambda}}\left(\widetilde{A} + \widetilde{B}K_e e + \widetilde{C}K_d \frac{d^\mu e}{dt^\mu}\right) + K_{PD}\left(\widetilde{A} + \widetilde{B}K_e e + \widetilde{C}K_d \frac{d^\mu e}{dt^\mu}\right) \\ &= \left[K_{PD}\widetilde{A} + K_{PI}\widetilde{A}\frac{t^\lambda}{\Gamma(\lambda+1)}\right] + \left[K_{PD}\widetilde{B}K_e\right]e + \left[K_{PD}\widetilde{C}K_d\right]\frac{d^\mu e}{dt^\mu} \\ &\quad + \left[K_{PI}\widetilde{B}K_e\right]\frac{d^{-\lambda}e}{dt^{-\lambda}} + \left[K_{PI}\widetilde{C}K_d\right]\frac{d^{\mu-\lambda}e}{dt^{\mu-\lambda}}\end{aligned} \tag{20}$$

Thus drawing an analogy with the classical PID controller structure, the first term $\left[K_{PD}\widetilde{A} + K_{PI}\widetilde{A}\dfrac{t^\lambda}{\Gamma(\lambda+1)}\right]$ represents a time dependent gain due to the presence of time in it. The term $\left[K_{PD}\widetilde{B}K_e\right]$ represents the proportional gain, $\left[K_{PD}\widetilde{C}K_d\right]$ represents the fractional order derivative gain, $\left[K_{PI}\widetilde{B}K_e\right]$ represents the fractional order integral gain and $\left[K_{PI}\widetilde{C}K_d\right]$ represents an additional FO integro-differential gain. The last term can represent either a fractional derivative or a fractional integral action depending on which value between $\{\lambda, \mu\}$ is greater.

Few recent research results show that band-limited implementation of FOPID controllers using higher order rational transfer function approximation of the integro-differential operators gives satisfactory performance in industrial automation e.g. in [39].

The Oustaloup's recursive approximation, which has been used to implement the integro-differential operator ($s^\alpha, \alpha \in \mathbb{R}$) in frequency domain is given by the following expression, representing a higher order analog filter.

$$s^\alpha \simeq K \prod_{k=-N}^{N} \frac{s + \omega'_k}{s + \omega_k} \quad (21)$$

where, the poles, zeros, and gain of the filter can be recursively evaluated using (22).

$$\omega_k = \omega_b \left(\frac{\omega_h}{\omega_b}\right)^{\frac{k+N+\frac{1}{2}(1+\alpha)}{2N+1}}, \omega'_k = \omega_b \left(\frac{\omega_h}{\omega_b}\right)^{\frac{k+N+\frac{1}{2}(1-\alpha)}{2N+1}}, K = \omega_h^\alpha \quad (22)$$

Thus, any signal $f(t)$ can be passed through the filter (21) and the output of the filter can be regarded as an approximation to the fractionally differentiated or integrated signal $D^\alpha f(t)$. In (21)-(22), $\alpha$ is the order of the differ-integration, $(2N+1)$ is the order of the filter and $(\omega_b, \omega_h)$ is the desired frequency range of fitting.

Even with the truncation of infinite dimensional natures of FO operators with high order IIR filters, the obtained FOPID controllers are found to outperform classical PID structure and even more with augmentation of fuzzy inferencing with the FO differ-integral operators [35], [40], [41]. Thus there is always a trade-off between the complexity of the realization of the FOPID controller and the achievable accuracy. In the present study, 5$^{th}$ order Oustaloup's recursive approximation is done for the integro-differential operators within a chosen frequency band of $\omega \in \{10^{-2}, 10^2\}$ rad/sec.

*3.2. Details of the fuzzy inference within the intelligent fractional order controller*

Fuzzy inference is the method by which the nonlinear mapping between the input and the output variables is established with the help of fuzzy logic. The process of fuzzy inferencing mainly comprises of fuzzy rule base, membership functions used in the rules, reasoning mechanism by the use of fuzzy logic operators, fuzzification and defuzzification operations etc. Here, the basic fuzzy logic controller uses a two dimensional rule base as shown in Figure 3 and triangular membership functions as reported in Das *et al.* [35] with 50% overlap.

| $\dfrac{d^\mu e}{dt^\mu}$ \ e | NL | NM | NS | ZR | PS | PM | PL |
|---|---|---|---|---|---|---|---|
| PL | ZR | PS | PM | PL | PL | PL | PL |
| PM | NS | ZR | PS | PM | PL | PL | PL |
| PS | NM | NS | ZR | PS | PM | PL | PL |
| ZR | NL | NM | NS | ZR | PS | PM | PL |
| NS | NL | NL | NM | NS | ZR | PS | PM |
| NM | NL | NL | NL | NM | NS | ZR | PS |
| NL | NL | NL | NL | NL | NM | NS | ZR |

**Figure 3: Rule base for error, error derivative and FLC output.**

Here, a Mamdani type inferencing is used with "min" type operator for implication and "max" type operator for rule aggregation. The error and its fractional derivative ($0 < \mu < 1$) is assumed to follow the rule base depicted in Figure 3 composed of 49 ($7 \times 7$) rules. The acronyms NL, NM, NS, ZR, PS, PM and PL refer to Negative Large, Negative Medium, Negative Small, Zero, Positive Small, Positive Medium and Positive Large respectively. The FLC outputs ($u_{FLC}$) in Figure 2 is derived with the center of gravity method for defuzzification. Fractional order enhancement of fuzzy logic based PID controller has been extensively studied by Das *et al.* [35]. The present paper applies the concept of tuning integro-differential orders along with the input-output SFs for the fuzzy FOPID structure using real coded genetic algorithm to ensure command following power level tracking of the nuclear reactor even with a shift in operating power level along with random sensor noise and stochastic delay consideration in the reactor control loop.

### *3.3. Objective function for optimization based controller tuning*

Comparison of different tuning strategies for $PI^\lambda D^\mu$ controllers are described in Das *et al.* [40], [41]. Among different approaches the time domain performance index optimization based tuning of fractional order controllers are quite easy and becomes essential when high nonlinearity, like the fuzzy inferencing in this case, comes into play and sophisticated frequency domain tuning techniques cannot be easily applied. The fuzzy $PI^\lambda D^\mu$ structure is tuned here with the minimization of a chosen control objective as a time domain performance index involving the control loop error and required controller effort. In the present study, the integral performance index ($J$) to be minimized has been taken as the weighted sum of Integral of Time multiplied Squared Error (ITSE) and Integral of Squared Controller Output (ISCO):

$$J = \int_0^\infty \left[ w_1 \cdot te^2(t) + w_2 \cdot u^2(t) \right] dt \tag{23}$$

Optimization result with objective function (23) produces the optimally tuned controller parameters (SFs and integro-differential orders) in terms of low error index and control signal. The inclusion of the squared error term in the ITSE penalizes the peak overshoot to a large extent. Also, the time multiplication term penalizes the error signal more at the later stages than at the beginning and hence results in a faster settling time. The squared controller output is also included in $J$ so that the control signal does not become too large and result in actuator saturation and integral windup. The weights $w_1$ and $w_2$ have been incorporated in the objective function (23) to keep a provision for balancing the impact of the error and the control signal. In this case the weights have been considered to be equal for the two objectives of (23) implying that the minimization of the error index and the control signal are equally important. A real coded GA has been employed here to tune the fuzzy FOPID controller while minimizing the custom objective function (23).

It is well known that genetic algorithm is a stochastic optimization process and is less susceptible to get trapped in local minima, compared to the gradient based algorithms. In GA, a solution vector is randomly chosen from the search space which undergoes reproduction, crossover and mutation in each iteration to give rise to a better population of solution vectors in the next iteration. The solution is refined iteratively until the objective function falls below a certain tolerance level or the maximum number of iterations are exceeded. Here, the number of population members in GA is chosen to be 20. The crossover fraction is taken to be 0.8 and the mutation fraction is 0.2. Also, the decision variables for GA are the input-output SFs and the integro-differential orders of the FO fuzzy PID controller i.e. $\{K_e, K_d, K_{PI}, K_{PD}, \lambda, \mu\}$.

## 4. Controller design at ideal condition with no consideration of stochastic phenomena in the reactor control loop

It is clear that the linearized reactor models (17) have widely varying dc gain which makes it difficult to design a single controller which can faithfully enforce command following power level tracking for different operating power of the reactor. In the first part of our exploration, the fuzzy $PI^\lambda D^\mu$ controller parameters have been tuned with GA for the reactor models corresponding to 100% and 20% of full power as in Table 1. The scaling factors of a conventional fuzzy PID controller with the integro-differential operators as unity ($\lambda = \mu = 1$) has also been tuned using the similar technique explored by Pan *et al.* [22] to get a fair comparison. Figure 4 shows the step responses of the nuclear reactor at different operating points with the optimum fuzzy $PI^\lambda D^\mu$ controller, tuned at the lowest and highest dc gain. Similar frequency domain fractional order controller tuning with lowest and highest power level models of nuclear reactor have been explored by Saha *et al.* [17] and Das *et al.* [18] respectively. It is evident from Figure 4 and Figure 5 that the set-point tracking performance is satisfactory for fuzzy $PI^\lambda D^\mu$ controller with smooth and small control signal which is desirable from an actuator design point of view. Whereas, with the conventional fuzzy PID controller, the time response becomes oscillatory with the 100% power controller and suffers from heavy oscillatory control signal with the 20% power controller that may cause damage to the actuator. Also it is evident that for both the controllers, the response is well-behaved when tuned with the full reactor power model since delayed tracking of the input reference for $PI^\lambda D^\mu$ controller can be observed from Figure 4. Similarly, it is not recommended to tune the conventional fuzzy PID controller at low (20%) power model since instantaneous jerks in the time response and oscillation in the control signal can be observed in Figure 5 even when there is no externally introduced stochastic consideration in the reactor control loop. Figure 5 clearly shows that the controller tuned at low reactor power performs inferior at higher power with fuzzy PID controller. The observation motivates to find the robust controller structure tuned at full power which can tolerate the change in reactor model due to nonlinearity i.e. reduction in operating power level when external disturbances are introduced in the control loop like network induced random delay and sensor noise. The simulation presented in this section clearly shows that the fuzzy $PI^\lambda D^\mu$ structure tuned with 100% reactor power model performs best amongst the four combinations i.e. two fuzzy controllers tuned at highest/lowest operating power.

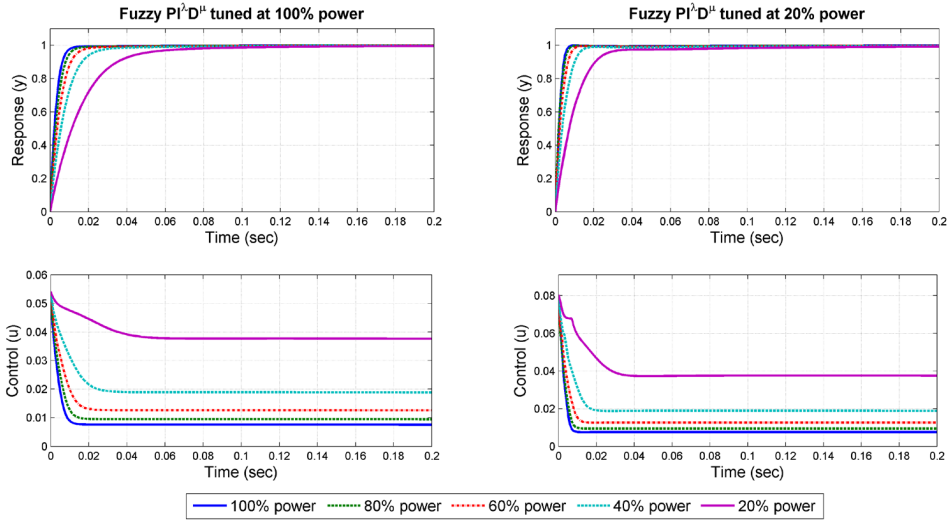

**Figure 4:** Performance with the Fuzzy $PI^\lambda D^\mu$ at different reactor power in ideal condition.

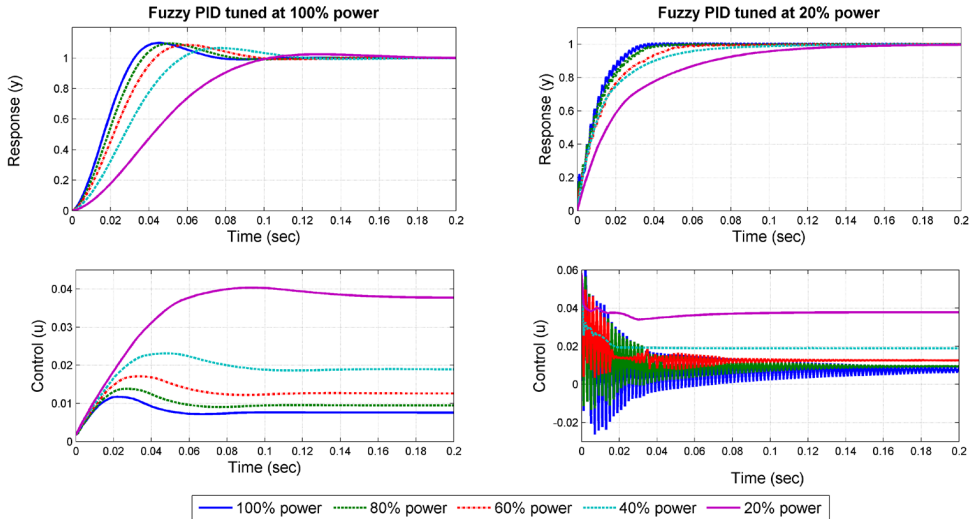

**Figure 5:** Performance with the Fuzzy PID at different reactor power in ideal condition.

## 5. Fuzzy FOPID controller design with stochastic phenomena in the control loop having long range dependence

### 5.1. Self-similar network induced delays in control loop

In this section, the assumption of ideal reactor operation is not made in order to see the relative capability of fuzzy PID and fuzzy $PI^\lambda D^\mu$ controller to suppress unwanted stochastic phenomena in the reactor control loop. For case study, a heavily loaded LAN data have been collected as shown in Figure 6 for estimation of the degree of self-similarity indicated by its associated Hurst parameter. Some typical characteristics of such network induced random delays in a loaded network have been discussed by Majumder *et al.* [42]. Detailed study regarding statistical characterization of delays in heavily loaded network has been done by Tejado *et al.* [30], [31]. From the time series of the random delays in a communication network it can be inferred that the random variable possesses some sudden spikes which may be confused with outliers and thus cannot be modeled using conventional Gaussian or mixture of Gaussian distributions. This typical case is known as non-Gaussian $\alpha$-stable distribution of the random variable which makes the mathematically tractable solutions for suppression of such stochastic phenomena with fractional order characteristics quite difficult [43]. Figure 7 shows that the run-time variance of such spiky random variables (packet delays) does not converge to a finite value and have a non-Gaussian $\alpha$-stable distribution [42].

The random network induced delay consideration has been motivated from the study reported in [42]. Majumder *et al.* [42] used various smoothing filters for removing the stochastic nature of the power signal which is continuously fed back from the reactor house to the RRS through a long communication medium. This contaminates the power signal measured by the Self-Powered Neutron Detectors (SPNDs) with such random delays as the instantaneous measurements are then sent to the controller through the long communication medium representing the feedback path of the reactor control loop. In this paper however, the objective is to design an efficient fuzzy logic based controller which can reject the unwanted stochastic phenomena in the control loop in a mean-squared sense and also ensure oscillation-free command following reactor maneuvering at all operating power levels, in spite of the presence of random delays and measurement noise in the sensors.

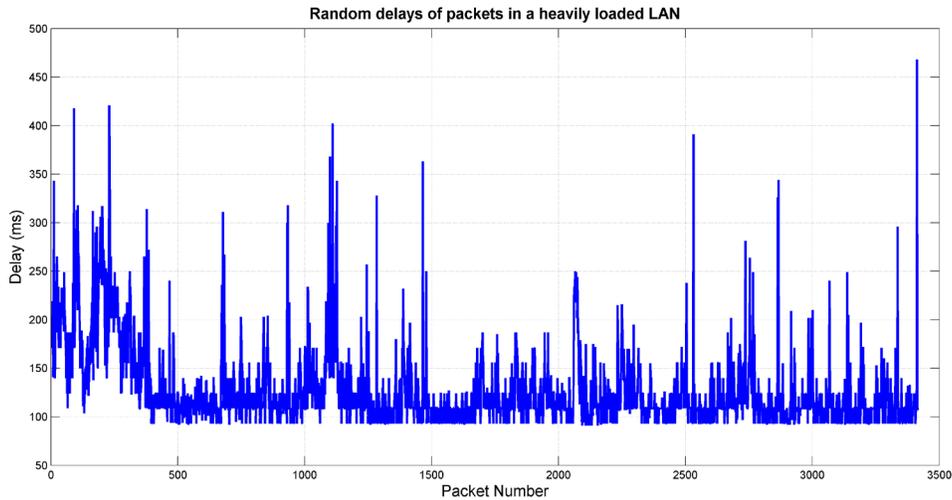

**Figure 6: Time domain presentation of the network induced stochastic delay.**

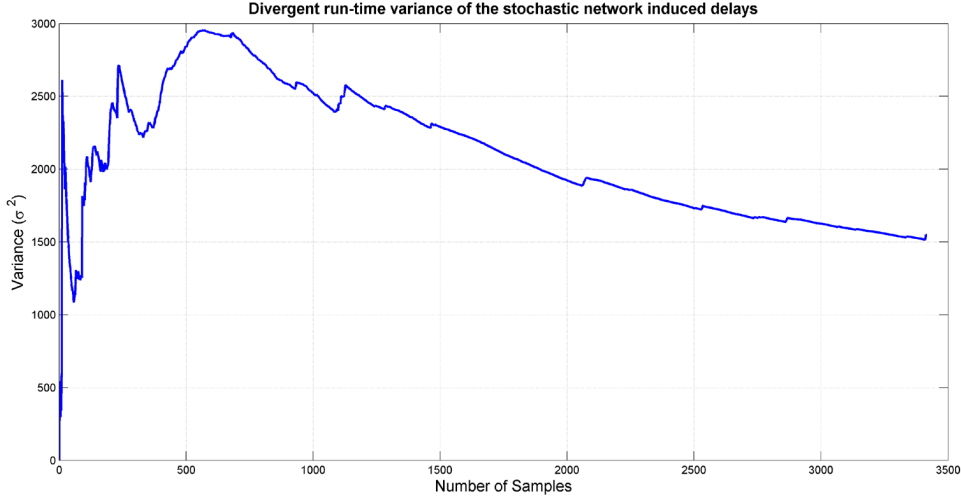

**Figure 7: Diverging run-time variance for the random network induced delay.**

## *5.2. Basics of long range dependency (LRD)*

The basic mathematical preliminaries of the random variables in the reactor control loop with long range interaction or dependence is discussed here. Let $X_t$ be a time series which is weakly stationary, suggesting that the series has a finite mean and the covariance depends only on the lag between the two points existing in the series. The self-similarity can be explained in different ways like using the concept of auto-correlation and power spectral density. Let $\rho(k)$ be the auto-correlation function (ACF) of $X_t$. The ACF $\rho(k)$ for a weakly stationary time series $X_t$ is given by $\rho(k) = E\left[\dfrac{(X_t - \mu)(X_{t+k} - \mu)}{\sigma^2}\right]$ with $E(X_t)$ being the statistical expectation of $X_t$. Also, $\mu$ and $\sigma^2$ are the mean and variance respectively. Self-similarity can best be related with long range dependency. The time series $X_t$ is said to have long range dependency if $\sum_{k=-\infty}^{k=\infty} \rho(k)$ diverges. Often $\rho(k)$ takes the form $\rho(k) \simeq C_\rho k^{-\alpha}$, with $C_\rho$ being positive and $\alpha \in (0,1)$. Parameter $\alpha$ is related to the Hurst parameter via the equation $\alpha = 2H - 1$. This is the most common definition of LRD. One very practical view-point of self-similarity can be revealed from the signal processing point of view i.e. its power spectral density. If the power spectral density of the time series $X_t$ be $f(\omega)$, then a time-series having LRD must conform to the following relation:

$$f(\omega) = \dfrac{\sigma^2}{2\pi} \sum_{k=-\infty}^{k=\infty} \rho(k) e^{ik\omega} \qquad (24)$$

where, $i = \sqrt{-1}$. This definition of spectral density comes from Wiener-Khintchine theorem. Also, the weakly stationary time series $X_t$ is said to have LRD if its power spectral density obeys $f(\omega) \simeq C_f |\omega|^{-\beta}$ with $C_f > 0$ and some real $\beta \in (0,1)$, where, $\beta$ is related to the associated Hurst parameter by the relation $H = (\beta+1)/2$.

*5.3. Hurst parameter and its estimation*

Many physical real world processes, exhibit LRD. Modeling of those physical processes require the correct estimation of the LRD which is measured by the Hurst parameter $H$. Estimation of LRD time series enhances the importance of analyzing the self-similarity in the time series. The concept of LRD was firstly introduced by Mandelbrot and Van Ness [44] in terms of Fractional Brownian Motion (FBM). Since then it has been addressed by many other contemporary researchers to analyze degree of self-similarity in a time series and leads to the concept of Hurst parameter. A second order time series $Y = f(u)$ is said to have a LRD if its auto-correlation function $\rho(\tau) = E[(f(\tau)f(0)]$ decays with the power law function of lag $\tau$ so that the auto-correlation series $\sum_{\tau} \rho(\tau)$ is not summable over the length of $\tau$. For the processes having Hurst parameter between $0 < H < 0.5$ are called anti-persistent process or negatively correlated. These processes generally have short range dependency. Processes with $0.5 < H < 1$ are called positively correlated. Hurst parameter $H = 0.5$ means the process is not correlated signifying conventional white-Gaussian noise. Processes with $1 < H < 1.5$ is said to have no dependency in time domain.

There are different methods to find out the Hurst parameter of a fractal time series. Popular method to find out the Hurst parameter is the R/S analysis. Apart from that there are a number of methods like aggregated variance method, absolute value method, Periodogram method, variance of residuals method, local whittle method, wavelet based method, Higuchi method and differenced variance approach etc. [43], [45]. Abrupt shift of mean in the time series or other contaminations can make the series non stationary and result in an over-estimation of Hurst parameter. In the present study, communication data taken from a heavily loaded shared network has been analyzed with one of the estimators as a test case. LRD can be thought of in two different ways. In time domain high degree of correlation in the distant samples of any time series can be modeled as LRD. In frequency domain significant amount of power at very low frequency can be an indication of presence of LRD. As reported in Majumder *et al.* [42], various estimators can be applied to find the Hurst parameter and the fractal dimension of the network induced random delays. In this paper the Rescaled range (R/S) analysis method has been chosen which has given the highest fractal dimension [42].

Let, $R(n)$ be the range of data aggregated over the block of length $n$ and $S(n)$ be the variance recorded over the same scale of range. For the series to follow self-similarity the following relation (25) must be maintained.

$$E\left[\frac{R}{S(n)}\right] \simeq C_H n^H \quad (25)$$

Taking logarithm on both sides of (25), the Hurst parameter can be estimated by the following regression formula (26).

$$\log E\left[\frac{R}{S(n)}\right] \simeq \log C_H + H \log n \tag{26}$$

From (26) it is clear that $C_H$ is a positive constant and independent of $n$. Thus, it should have a constant slope as $n$ becomes large. If the process sample is drawn from a stable distribution, slope of (26) will have the value of 0.5 like the random Gaussian noise. If the slope is over 0.5 it indicates the persistency in the time series. If the slope is below 0.5 an ergodic mean reverting process is indicated. Small value of $n$ will make the result anti-persistent and the resulted value of Hurst parameter will be invalid. Again a large value of $n$ will produce too few samples to correctly estimate $H$ for a self-similar process. Therefore, choice of $n$ should be judicious with this type of estimator. For the network induced delay shown in Figure 6, the R/S analysis gives the indication of long range dependency as Hurst $H = 0.8837$, fractional order $\alpha = 2H - 1 = 0.7674$, fractal dimension $D = 2 - H = 1.2326$. Here, the fractional order of the random delay implies a $1/f^\alpha$ process having self-similar nature.

*5.4. Persistent and anti-persistent noise*

The concept of persistent and anti-persistent noise is somewhat similar to the above mentioned preliminaries regarding Hurst parameter of networked delay. But in order to clearly distinguish between random network-induced delay and noise both having long range dependence, the generation of fractional Gaussian noise has been illustrated here with slightly different mathematical notations.

The stationary sequence $x(nT)$ obtained by sampling the fractional Brownian motion process $B_H(t)$ with a sampling interval $T$ and then calculating the first difference is known as the discrete fractional Gaussian noise. It can be represented as

$$x(nT) \triangleq B_H(nT) - B_H(nT - T) \tag{27}$$

As the statistical properties of the fractional Brownian motion (fBm) do not change with scale, we set $T = 1$ to obtain the discrete fractional Gaussian noise process as $x(n) \triangleq B_H(n) - B_H(n-1)$ and it is referred to as fGn. The autocorrelation of the discrete fractional Gaussian noise is

$$r_x(l) = \frac{1}{2}\sigma_H^2 \left(|l-1|^{2H} - 2|l|^{2H} + |l+1|^{2H}\right) \tag{28}$$

The fGn is wide-sense stationary as the correlation depends only on the distance $l$ between the samples. Also for $H = 1/2$, we have $r_x(l) = \delta(l)$ which implies that the fGn process is white noise.

The fGn process has long memory in the range $1/2 < H < 1$, since $\sum_{l=-\infty}^{\infty} r(l) = \infty$ or equivalently $R(e^{j\omega}) \to \infty$ as $|\omega| \to 0$. In this case the autocorrelation decays slowly and the frequency response is analogous to a low-pass filter. However, the process exhibits

short memory for $0 < H < 1/2$, since $\sum_{l=-\infty}^{\infty} |r(l)| < \infty$ and $\sum_{l=-\infty}^{\infty} r(l) = 0$, or equivalently $R(e^{j\omega}) \to 0$ as $|\omega| \to 0$. Also, for $0 < H < 1/2$, the correlation is negative, that is, $r(l) < 0$ for $l \neq 0$ and hence the process exhibits anti-persistence. For this instance, the autocorrelation decays very quickly and the frequency response is analogous to that of a high-pass filter.

The fGn can alternatively be viewed as the output of a fractional integrator ($1/s^{\beta}, 0 < \beta < 1$) driven by continuous time stationary white Gaussian noise (wGn) $w(t)$ with variance $\sigma^2$. The output can be written in the following convolution form:

$$y_{\beta}(t) = \frac{1}{\Gamma(\beta)} \int_{-\infty}^{t} w(\tau)(t-\tau)^{\beta-1} d\tau \tag{29}$$

There are many methods of approximately simulating the self-similar nature of the noise on a finite bit computer like the spectral synthesis method, multiple time scale fluctuation approach, fast fractional Gaussian noise generator etc. [43], [45]. As is evident this would only be an approximation limited by the precision of the computations as self-similarity would not exist at all scales. In this paper the fGn is generated by passing a wGn it through a fractional order integrator as in equation (21).

From the study by Chen *et al.* [46], it is seen that in electrochemical processes the Hurst parameter of the noise has been found out to be $H = 0.834$ or in other word the fractional order of the noise is $\beta = 2H - 1 = 0.668$, where the noise can be visualized as a $1/f^{\beta}$ process. This clearly indicates that in many naturally occurring chemical processes the realistic random noise is a filtered version somewhere between wGn ($\beta = 0$) and pink noise ($\beta = 1$) which can be obtained by passing it through an integrator. We generated fGn using the Oustaloup's filter (21) representation of FO integrator similar to that used in the fuzzy FOPID controller simulations. In the present study, the SPND signal of nuclear reactor power level is assumed to be contaminated by measurement noise having a noise with $\beta = 0.668$ denoting its persistent nature. For the consideration of anti-persistent noise, the sign of the fractional order is just made opposite denoting a negative correlation of the noise samples. Similar studies have been done by Pan *et al.* [33] to suppress persistent and anti-persistent noise using fractional order fuzzy PID controller which is also enhanced in the present work with the introduction of network induced delay as well and taking the shift in operating power level of the nuclear reactor into consideration.

### *5.5. Performance degradation with the controllers designed at ideal condition*

The fuzzy PID and the FO Fuzzy PID controllers are tuned at 100% and 20% reactor power without the consideration of the network induced delays and the sensor noise in section 4. Now the well-tuned control loops are tested at other operating conditions as well with the stochastic network delays and the noise present and also having long range dependence given by the statistics reported in section 5.1-5.4 for the respective cases. Figure 8-Figure 10 show the set point tracking and control signal of the fuzzy PID controllers under these circumstances. A few important conclusions can be drawn from the figures. The fuzzy PID controller tuned at 100% power has a good set point tracking and control signal at the other operating conditions with lower power. However if the same fuzzy PID controller is tuned at a 20% power then it shows poor performance as the power level increases. Especially at 100% power, the controller shows highly noisy response in both the output and the control signal. Thus the controller tuning should always be at the maximum power to ensure sufficient robustness when the nuclear reactor is operating at other conditions as well. A comparison of Figure 8-Figure 10 also shows that the effect of anti-persistent noise is much more detrimental than the persistent or the white Gaussian noise. The process output of the fuzzy PID controller tuned at 100% power is also affected significantly which is not observed when the persistent noise and the white Gaussian noise are present. However, more than the process output i.e. the power level, the control signal or the manipulated variable i.e. the change in reactivity due to control rod movement is highly affected due to the anti-persistent nature of the noise. Random fluctuations in the control signal or the position of control rods are highly detrimental and might result in actuator failure, mechanical shocks to the actuation system and tripping of the nuclear power plant altogether.

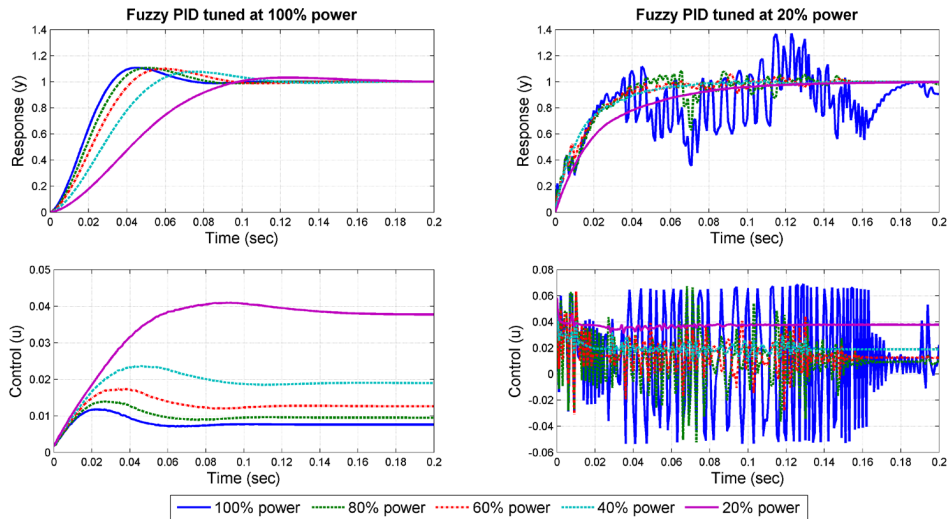

**Figure 8: Nominal tuning of Fuzzy PID controller to handle self-similar network delay (α=0.7674) and persistent sensor noise (β=0.668).**

Figure 11-Figure 13 show the FO fuzzy PID controllers tuned at 100% power and 20% power without considering the random network delay or the noise in the control loop. It can be seen that the Fuzzy FOPID controllers work well for all operating conditions irrespective of whether they are tuned at 100% power or at 20% power level, unlike that of the fuzzy PID controllers. In the presence of network delays and persistent noise in Figure 11 both the set point tracking and the control signal is good without any jittery effects. In the presence of network delays and white Gaussian noise in Figure 12, there are small jittery effects for both the controllers (i.e. one tuned at 100% and one at 20%) and this jittery effect is more pronounced in Figure 13 when there is anti-persistent noise along with network delays. A comparison with the curves obtained for the Fuzzy PID controller in Figure 8-Figure 10, show that the fuzzy FOPID is capable of suppressing oscillations of the power level and the control variables much more than the simple fuzzy PID controller, hence showing better robustness characteristics.

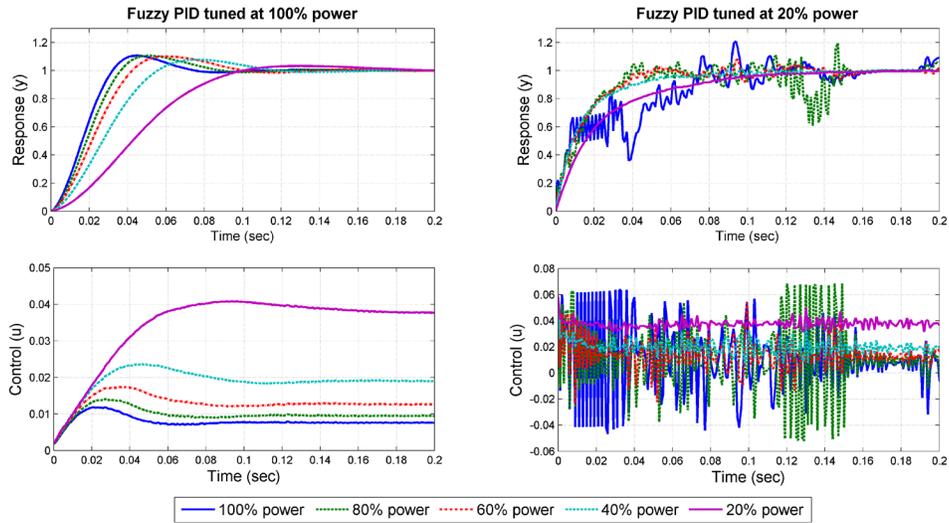

**Figure 9: Nominal tuning of Fuzzy PID controller to handle self-similar network delay (α=0.7674) and white Gaussian sensor noise (β=0).**

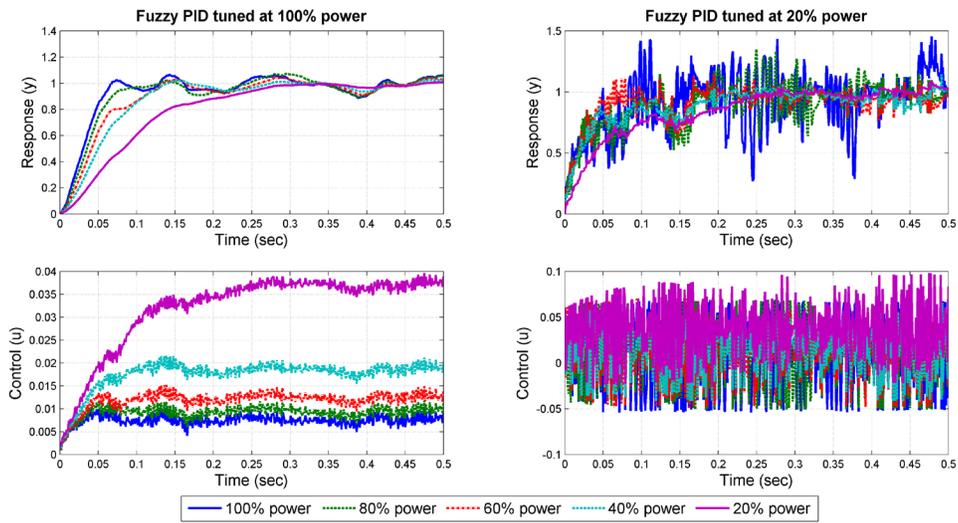

**Figure 10: Nominal tuning of Fuzzy PID controller to handle self-similar network delay (α=0.7674) and anti-persistent sensor noise (β=-0.668).**

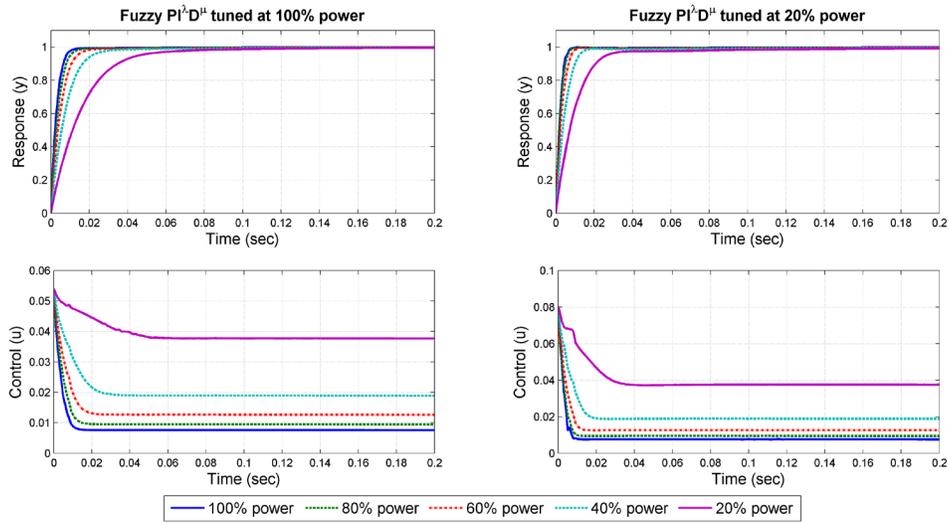

**Figure 11: Nominal tuning of Fuzzy FOPID controller to handle self-similar network delay (α=0.7674) and persistent sensor noise (β=0.668).**

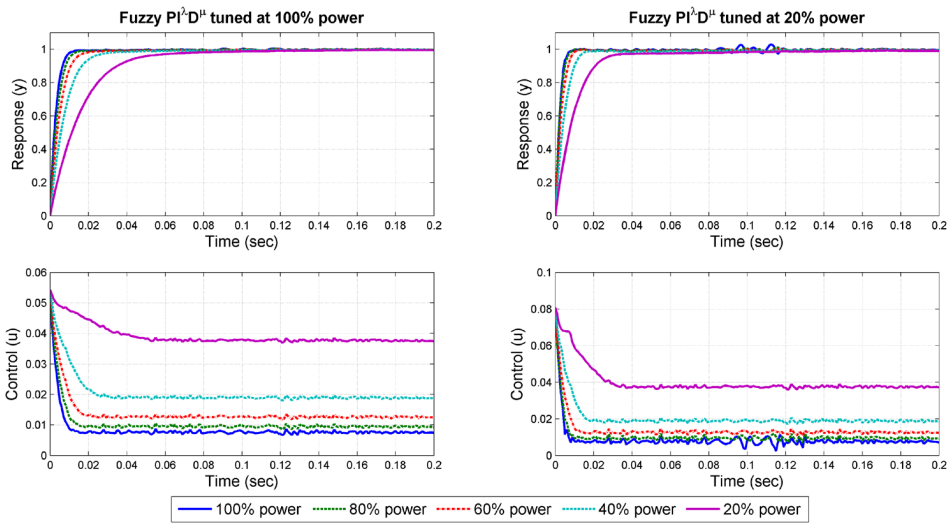

**Figure 12: Nominal tuning of Fuzzy FOPID controller to handle self-similar network delay (α=0.7674) and white Gaussian sensor noise (β=0).**

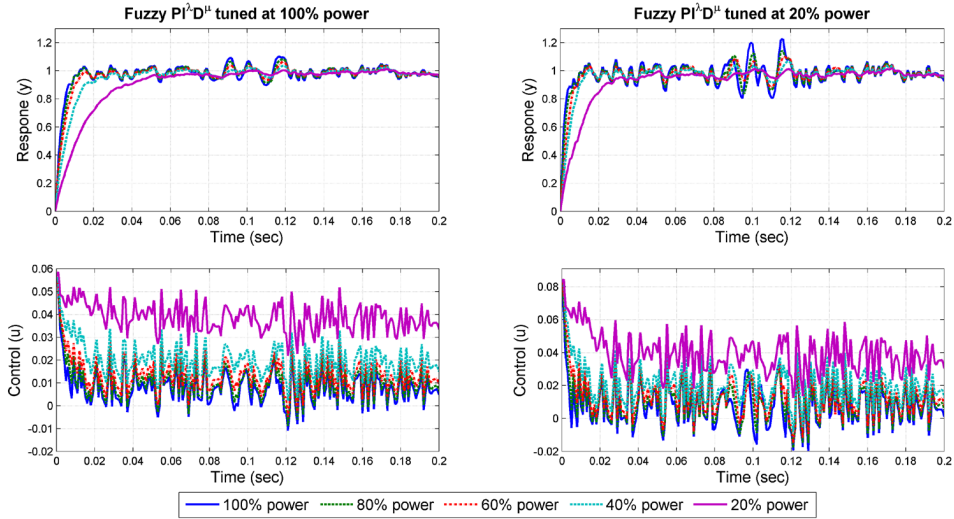

**Figure 13: Nominal tuning of Fuzzy FOPID controller to handle self-similar network delay (α=0.7674) and anti-persistent sensor noise (β=-0.668).**

## 6. Improved stochastic optimization based controller design with random noise and delay consideration

### 6.1. *Advantage of improved stochastic tuning of fuzzy FOPID controller*

The simulations as reported in section 5.5 are redone again with an improved tuning methodology as reported in Pan *et al.* [25], [33] for both the simple and the FO fuzzy PID controller. The GA based optimization of fuzzy controller SFs are done considering the stochastic variation of the objective function due to network delay and sensor noise. The concept is somewhat similar to the minimization of an expected value of the objective function using GA. It has been shown by contemporary researchers [47]-[49] that evolutionary and swarm based algorithms naturally works well for such cases where the objective function is time varying and noisy. The corresponding tuning results for the fuzzy PID and FOPID controllers are reported in Table 2 for different operating power level and stochastic considerations due the sensor noise or network delay.

**Table 2: Fuzzy controller Parameters at 100% & 20% Power in ideal condition (no sensor noise or network delay)**

| Power (%) | Controller | $J_{min}$ | $K_e$ | $K_d$ | $K_{PI}$ | $K_{PD}$ | $\lambda$ | $\mu$ |
|---|---|---|---|---|---|---|---|---|
| 100 | $PI^\lambda D^\mu$ | 0.0002 | 0.3236 | 0.0683 | 4.557 | 0.124 | 0.9643 | 0.0958 |
| 100 | PID | 0.0003 | 0.9918 | 0.0061 | 1.2510 | 0.001 | - | - |
| 20 | $PI^\lambda D^\mu$ | 0.0054 | 0.6534 | 0.3349 | 2.189 | 0.092 | 0.8407 | 0.0254 |
| 20 | PID | 0.0055 | 0.9859 | 0.0059 | 1.876 | 0.067 | - | - |

While evaluating the objective functions with stochastic phenomena, the control system is simulated with self-similar network induced delay and persistent/anti-persistent noise. This makes the objective function random as has been explicitly illustrated in [25], [33]. Thus for the same controller gains, the objective function will take different values due to the stochastic nature of the network delays and the noise. Hence, for each controller gain, the objective function is evaluated multiple times and the expected value is taken by the optimization algorithm. Figure 14-Figure 16 show the performance of the fuzzy PID controller. Similar characteristics can be seen for the improved tuning methodology as well. The tuned controllers work the best in the presence of the persistent noise and the network induced delays. The anti-persistent noise is the most detrimental for the controller's performance. However the most significant difference is that the band of oscillations is much reduced using the improved tuning methodology over the simple methodology reported in section 4 where the latter does not take the stochastic phenomena into consideration.

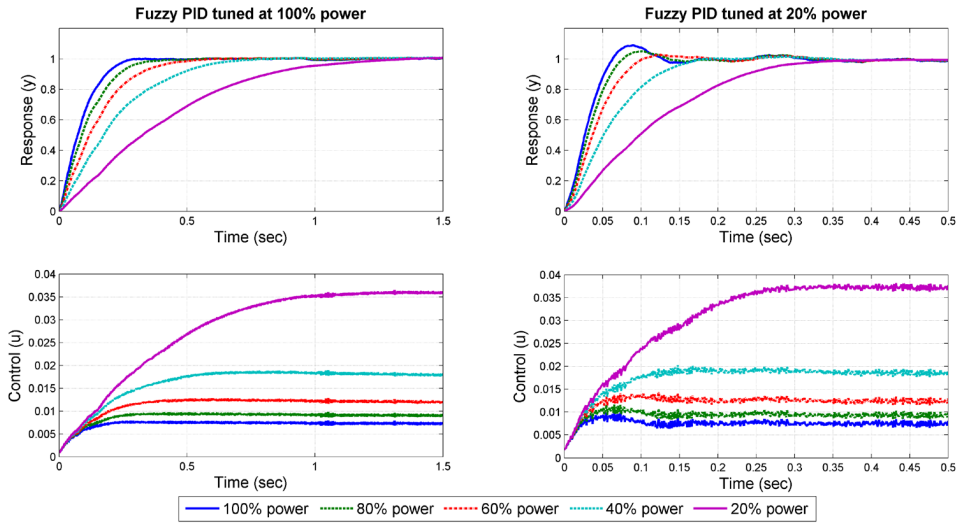

**Figure 14:** Improved tuning of Fuzzy PID controller to handle self-similar network delay ($\alpha$=0.7674) and persistent sensor noise ($\beta$=0.668).

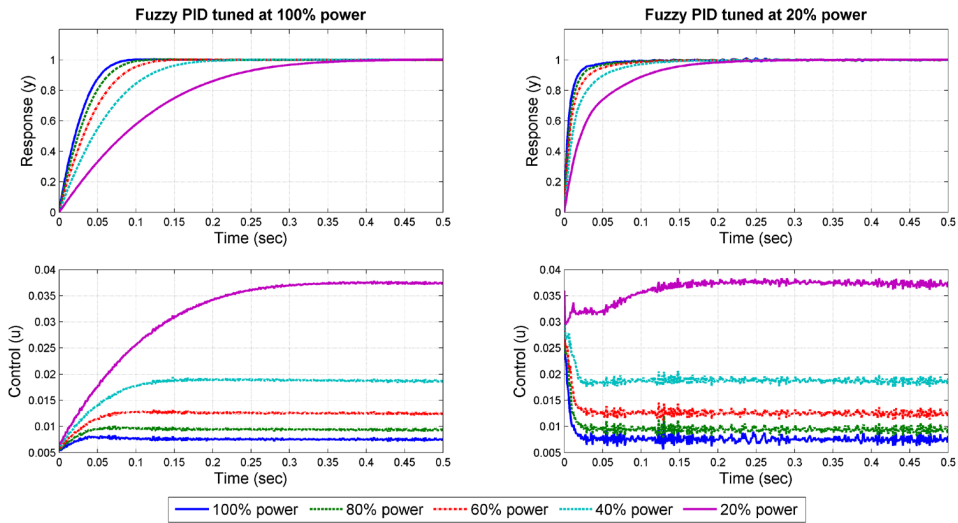

**Figure 15:** Improved tuning of Fuzzy PID controller to handle self-similar network delay ($\alpha$=0.7674) and white Gaussian sensor noise ($\beta$=0).

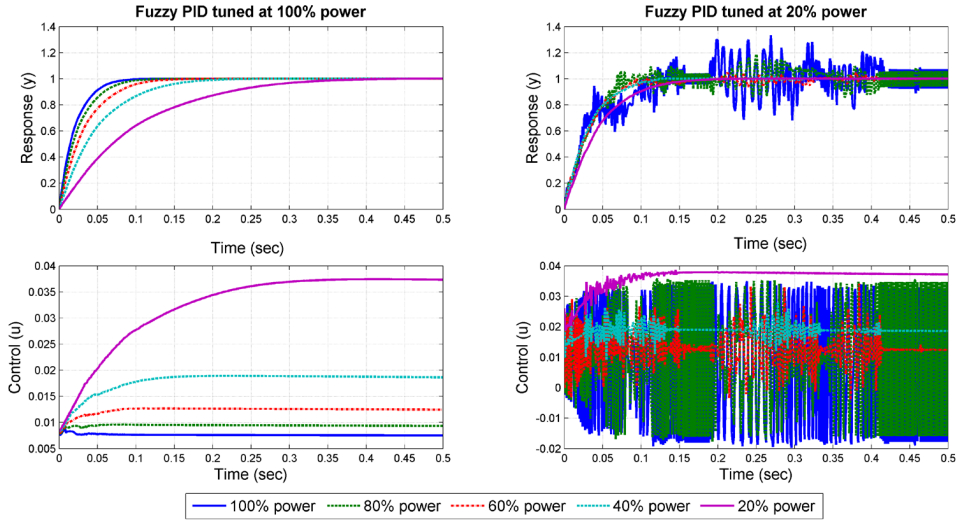

**Figure 16: Improved tuning of Fuzzy PID controller to handle self-similar network delay (α=0.7674) and anti-persistent sensor noise (β=-0.668).**

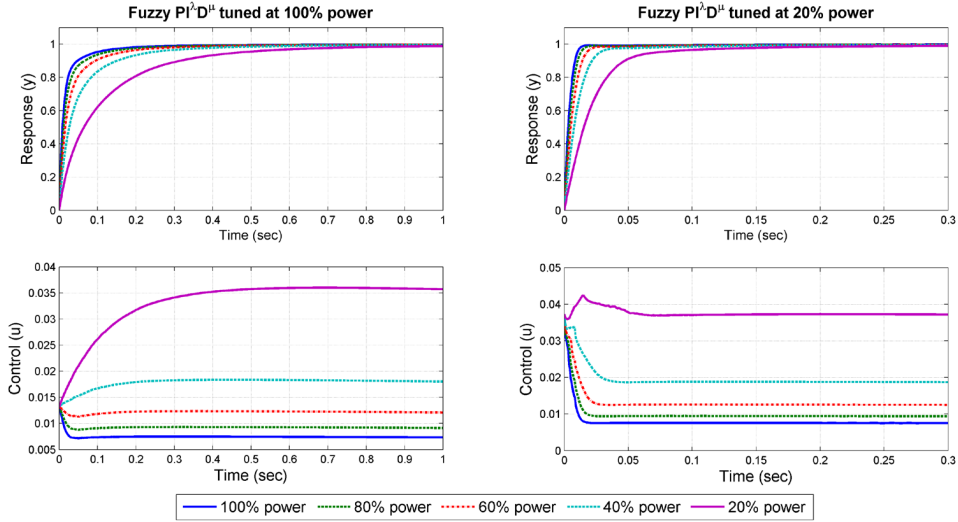

**Figure 17: Improved tuning of Fuzzy FOPID controller to handle self-similar network delay (α=0.7674) and persistent sensor noise (β=0.668).**

Figure 17-Figure 19 show the performance of the FO fuzzy controllers handling self-similar network delays and noise with long range dependence. As can be seen from the figures, the noise suppression is drastically improved for the fuzzy FOPID controllers over the simple fuzzy PID controllers with this improved tuning method. The pictorial representations of the improvement are quantitatively backed by the numerical results in Table 3. In Figure 17, the FO fuzzy controller shows almost no oscillations in the set point tracking or the control signal in the presence of persistent noise unlike the other previous simulations. The simulation with white Gaussian noise and network delay in Figure 18 shows slight oscillations in the control signal but the set-point tracking has almost no induced oscillations in it. The anti-persistent noise is the most detrimental as in the previous simulations as well and there are small oscillations in the control signal as well as the reference signal. However these oscillations are significantly smaller as compared to the previous cases, thus proving the merit of the improved tuning methodology and the advantage of fuzzy FOPID controller.

**Table 3: Fuzzy PID and $PI^\lambda D^\mu$ Controller parameters tuned with stochastic consideration of the reactor control loop**

| Stochastic consideration | Power (%) | Controller structure | $J_{min}$ | $K_e$ | $K_d$ | $K_{PI}$ | $K_{PD}$ | $\lambda$ | $\mu$ |
|---|---|---|---|---|---|---|---|---|---|
| Persistent noise ($\beta$=0.668) and self-similar delay ($\alpha$=0.7674) | 100 | Fuzzy $PI^\lambda D^\mu$ | 0.0003 | 0.1145 | 0.0796 | 0.9951 | 0.0488 | 0.9361 | 0.0863 |
| | | Fuzzy PID | 0.0511 | 0.9624 | 0.0011 | 0.1245 | 0.0011 | - | - |
| | 20 | Fuzzy $PI^\lambda D^\mu$ | 0.0055 | 0.6436 | 0.1946 | 1.4269 | 0.0514 | 0.8779 | 0.0569 |
| | | Fuzzy PID | 0.0561 | 0.9708 | 0.0021 | 0.5439 | 0.0019 | - | - |
| White noise ($\beta$=0) and self-similar delay ($\alpha$=0.7674) | 100 | Fuzzy $PI^\lambda D^\mu$ | 0.0005 | 0.1674 | 0.0010 | 1.6514 | 0.0361 | 0.9408 | 0.7431 |
| | | Fuzzy PID | 0.0006 | 0.3229 | 0.0010 | 1.1039 | 0.0212 | - | - |
| | 20 | Fuzzy $PI^\lambda D^\mu$ | 0.0058 | 0.5916 | 0.1182 | 1.0812 | 0.0580 | 0.8823 | 0.0010 |
| | | Fuzzy PID | 0.0058 | 0.9979 | 0.0010 | 0.9718 | 0.0402 | - | - |
| Anti-persistent noise ($\beta$=-0.668) and self-similar delay ($\alpha$=0.7674) | 100 | Fuzzy $PI^\lambda D^\mu$ | 0.0546 | 0.0802 | 0.4775 | 0.1638 | 0.0075 | 0.9866 | 0.0148 |
| | | Fuzzy PID | 0.0003 | 0.9166 | 0.0010 | 0.4412 | 0.0109 | - | - |
| | 20 | Fuzzy $PI^\lambda D^\mu$ | 0.0562 | 0.0928 | 0.2729 | 0.9870 | 0.0122 | 0.9875 | 0.0036 |
| | | Fuzzy PID | 0.0059 | 0.9985 | 0.0166 | 1.5201 | 0.0291 | - | - |

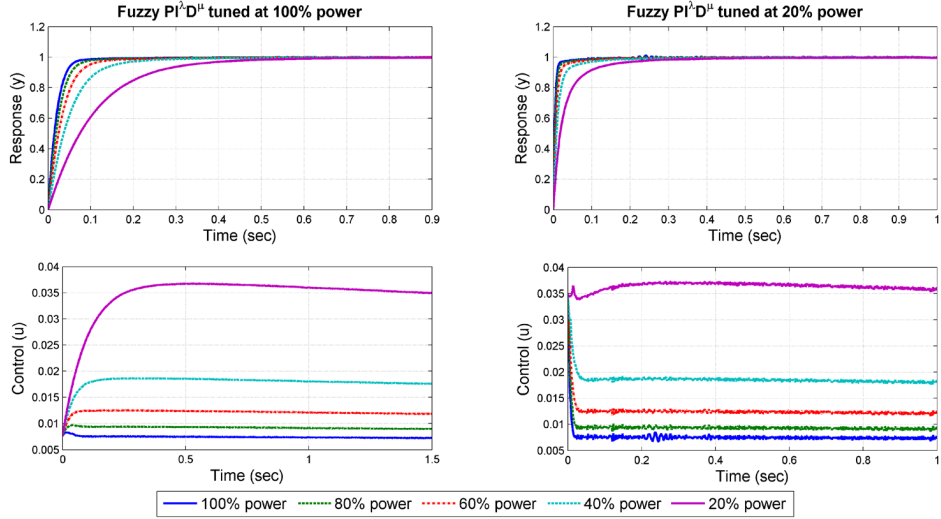

**Figure 18: Improved tuning of Fuzzy FOPID controller to handle self-similar network delay ($\alpha$=0.7674) and white Gaussian sensor noise ($\beta$=0).**

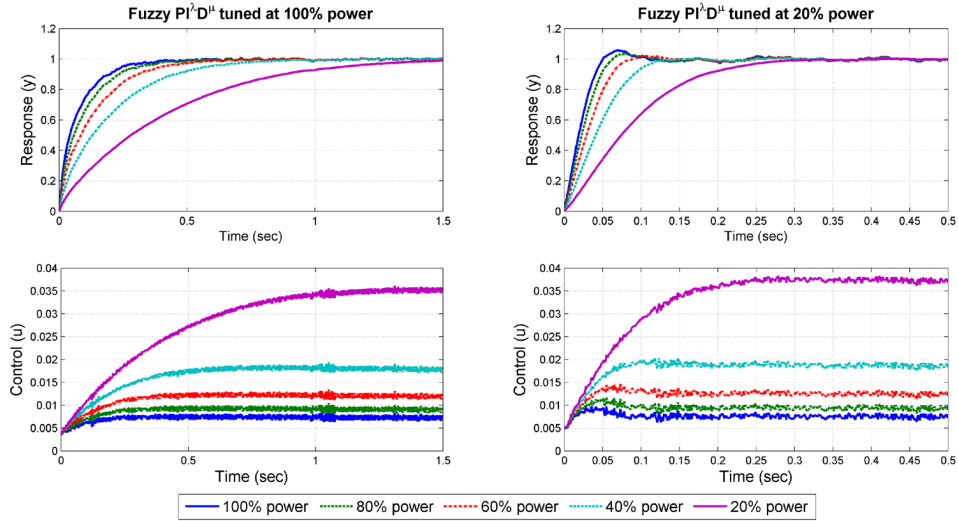

**Figure 19:** Improved tuning of Fuzzy FOPID controller to handle self-similar network delay ($\alpha=0.7674$) and anti-persistent sensor noise ($\beta=-0.668$).

### *6.2. Superiority of tuning fuzzy FOPID controller using the reactor model at full power*

From the above mentioned simulations it can concluded that the fuzzy controllers tuned at full power gives better robustness at other operating points with or without the random delay and noise in the control loop compared to the controllers tuned with reactor models at part load. This is a typical improvement over the technique adopted by Saha *et al.* [17] and in harmony with the finding of Das *et al.* [18]. Also, Pan *et al.* [25] have shown that such stochastic optimization based fractional order controller design does not affect the system's time response in a great extent with variation in the distribution of random phenomena in the control loop. But present simulation shows that it does, depending on its associated Hurst parameter or degree of long range dependence which is also in harmony with the findings in [33]. This is due the fact that the variance of the randomness seems to get increased in a significant extent for anti-persistent noise and it decreases for persistent noise, compared to the variation in distribution of the network induced stochastic phenomena like delay and packet drops as studied in [25]. Also the fractional order fuzzy PID controller tuned at 100% reactor power and the worst case of randomness in the loop i.e. anti-persistent noise and self-similar delay has the highest robustness. Since it is tuned for the worst case scenario, the fuzzy FOPID automatically handles rest of the relatively low complex cases like the ideal case with no network delay and noise, only self-similar network delay, self-similar network delay and white Gaussian noise, self-similar network delay and persistent noise etc. As discussed above, the random self-similar network delay and noise should not be confused to be same since the former introduces no extra energy in the control loop. Even though these two stochastic phenomena in reactor control loop has highly detrimental effect but the proposed stochastic tuning of fuzzy FOPID controller is faithfully capable of suppressing such unwanted effects at a wide range of reactor operation. Such self-similar consideration of the sensor noise must be done while designing controller for the RRS loop, as the early studies by Mandelbrot and Van Ness [44] have shown that fission reaction has intrinsically 1/f nature. This may contaminate the feedback signal sensed by SPND and must be filtered using some advanced control or signal processing technique to avoid any possible failure of safety critical applications like control systems for maneuvering nuclear reactor power.

### 6.3. Justification of choosing complex FO fuzzy control scheme for nuclear power level monitoring

It might seem that the adopted control structure is quite complex compared to the system under control i.e. the PWR models. The nonlinear state-space model of the reactor may seem to be simplistic but it should be noted that its linearized version at different operating powers are quite difficult to handle with a single controller. Since it is not feasible to switch the controller every time or tune it several times, depending on the shift operating point, a robust control strategy is required which will enforce guaranteed power level tracking at all operating points. Designing robust controller for nuclear reactor power level monitoring, while considering models at different operating point is not a trivial problem and even so when stochastic disturbances in the form of noise and network delays are present in the control loop. Please see investigations in [17] and [18] in related topics dealing with frequency domain robust controller design for different operating points of nuclear reactor.

Also, from the transfer functions in equation (17), it is evident that the linearized models have widely varying dc gain, co-existence of very fast and very slow time constants and strong leads. At all of the operating powers, one pole is very close to the origin which will make the open loop system behave somewhat similar to a marginally stable system [36], though the original system was governed by a set of nonlinear differential equations. Tiwari *et al.* [2] have shown that for large reactors even the open loop system may be unstable due to the presence of right half plane poles. Das *et al.* [35] gave a detailed simulation comparison amongst PID, FOPID, fuzzy PID and fuzzy FOPID controllers to handle a nonlinear plant and an open loop unstable plant, both of which types are present in the model of a nuclear reactor. Since, in [35] it has been already shown that the complex controllers especially having fuzzy inferencing give good performance for such complex systems, here we have only shown the comparison of fuzzy PID and fuzzy FOPID structures. Also, the effectiveness of fuzzy PID [22] and fractional order $PI^\lambda D^\mu$ controller [25] to suppress network induced stochastic disturbances, over simple control structures has made their hybrid i.e. fuzzy FOPID control scheme an obvious choice.

It is also understandable that many process control loops e.g. the light water reactors are controlled using PI type controllers since they do not amplify noise. Addition of derivative controller though increases the sensitivity against noise but helps in improving the stability of the open loop integrating/unstable reactor model by addition of zeros. Since the controller used here is not of conventional PID type and most popular fuzzy controllers needs the information about the rate of change of the error for deriving the control law by fuzzy inferencing, the derivative action in the controller structure cannot be avoided here. But the amplification of the stochastic phenomena can be reduced by using fractional derivative of error, smoothing effect of integer and fractional integrator after the FLC and averaging nature of the FLC.

### 7. Conclusion

Fractional order fuzzy PID controllers have been tuned for optimal power tracking of a nuclear reactor in load following mode. The reactor model is developed from the basic point-kinetics with reactivity feedback and thermal-hydraulic considerations. GA based optimization is used to tune the fuzzy $PI^\lambda D^\mu$ controller parameters to achieve efficient control performance and is illustrated by credible simulations. The fractional order fuzzy PID controller is shown to be robust when tuned at any power level, whereas the simple fuzzy PID controller must be tuned at the highest power level to show sufficient robustness at other operating conditions. Simulation results also show that among the different type of noise, anti-persistent noise in the control loop is the most detrimental from the performance point of view. The superiority of the fractional order fuzzy PID controller over the simple fuzzy PID controller in terms of noise suppression and robust set-point tracking is demonstrated by simulation examples. An improved controller tuning method is proposed using stochastic optimization which takes the random network delay and the noise into consideration in the tuning phase itself. Simulation results show improved performance of both the simple fuzzy PID and the fuzzy FOPID controller using the proposed methodology. Future scope of research can be directed towards the analysis of power level regulation in large nuclear reactors with spatial oscillations as well.